# PREHOMOGENEOUS VECTOR SPACES AND FIELD EXTENSIONS III

Akihiko Yukie[1]

Oklahoma State University

**Introduction**

Throughout this paper, $k$ is a field of characteristic zero, and $\overline{k}$ is its algebraic closure. We first recall the definition of prehomogeneous vector spaces.

**Definition (0.1)** *Let $G$ be a connected reductive group, $V$ a representation of $G$, and $\chi$ a non-trivial character of $G$, all defined $k$. Then $(G, V, \chi)$ is called a prehomogeneous vector space if it satisfies the following properties.*
*(1) There exists a Zariski open orbit.*
*(2) There exists a polynomial $\Delta(x) \in k[V]$ such that $\Delta(gx) = \chi(g)^a \Delta(x)$ for a certain positive integer $a$.*

Such $\Delta$ is called a relative invariant polynomial. We define $V^{\mathrm{ss}} = \{x \in V \mid \Delta(x) \neq 0\}$ and call it the set of semi-stable points. If $(G, V, \chi)$ is an irreducible representation, the choice of $\chi$ is essentially unique and we may write $(G, V)$ as well. For $x \in V_k^{\mathrm{ss}}$, let $G_x$ be the stabilizer and $G_x^0$ its connected component of 1.

The orbit space $G_k \setminus V_k^{\mathrm{ss}}$ for prehomogeneous vector spaces usually parametrizes interesting arithmetic objects and has been considered by many people. For example Gauss [1] considered the space of binary quadratic forms to investigate ideal classes of quadratic fields. Igusa [3] investigated spinors of dimension up to twelve, and these cases turned out to be special cases of prehomogeneous vector spaces. Later Igusa [4] formulated the problem in terms of Galois cohomology, and this formulation was used to parametrize field extensions of degree up to five in [11].

Let $W = k^7$ be the standard representation of $\mathrm{GL}(7)$ (i.e., the space of seven dimensional column vectors), and $e_1, \cdots, e_7$ the standard coordinate vectors. We use the notation $e_{ijk} = e_i \wedge e_j \wedge e_k$, etc. We consider the natural representation of $\mathrm{GL}(7)$ on $\wedge^3 W$. Let

$$\overline{w} = e_{234} + e_{567} + e_1 \wedge (e_{25} + e_{36} + e_{47}) \in \wedge^3 W.$$

We define $G_1$ to be the set of $g \in \mathrm{GL}(7)$ such that $g\overline{w} = c(g)\overline{w}$ for some $c(g) \in \mathrm{GL}(1)$. Then $G_1$ becomes a closed subgroup of $\mathrm{GL}(7)$, and $c : G_1 \to \mathrm{GL}(1)$ is a character of $G$ both defined over $k$.

In this paper, we consider the following two prehomogeneous vector spaces

---

[1]The author is partially supported by NSF grant DMS-9401391

(1) $G = G_1$, $V = W$,
(2) $G = G_1 \times \mathrm{GL}(2)$, $V = W \otimes k^2$.

**Definition (0.2)** $\mathfrak{E}\mathfrak{x}_i$ *is the set of isomorphism classes of Galois extensions of $k$ which are splitting fields of degree $i$ equations.*

In §1, we recall the relation between the orbit space $G_k \setminus V_k^{\mathrm{ss}}$ and Galois cohomology. In §§2,3, we construct a natural map $\alpha_V : G_k \setminus V_k^{\mathrm{ss}} \to \mathfrak{E}\mathfrak{x}_2$ for cases (1), (2) and prove that $\alpha_V$ is surjective. For case (2), let $\widetilde{T} = \mathrm{Ker}(G \to \mathrm{GL}(V))$. In §4 (Theorems (4.4), (4.13)), we determine the fiber structure of $\alpha_V$ and prove the following theorem.

**Theorem (0.3)** (1) *For case (1), the map $\alpha_V$ is bijective. Moreover, if $x \in V_k^{\mathrm{ss}}$ corresponds to a quadratic extension $k(\alpha)/k$ with $\alpha^2 \in k$, $G_x^0 \cong \mathrm{SU}_\alpha(2,1)$.*
(2) *For case (2), $\alpha_V^{-1}(k)$ consists of a single orbit, and if $k(\alpha)/k$ is a quadratic extension with $\alpha^2 \in k$, $\alpha_V^{-1}(k(\alpha)) \cong (\mathbb{Z}/2\mathbb{Z}) \setminus k^\times / \mathrm{N}_{k(\alpha)/k}(k(\alpha)^\times)$ where $\mathbb{Z}/2\mathbb{Z}$ acts on $k^\times / \mathrm{N}_{k(\alpha)/k}(k(\alpha)^\times)$ by $s \to s^{-1}$. Moreover, if $x \in \alpha_V^{-1}(k(\alpha))$ corresponds to $s \in k^\times$, $G_x^0/\widetilde{T}$ is the unitary group associated with the Hermitian matrix $\begin{pmatrix} -s & \\ & 1 \end{pmatrix}$.*

For the definition of $\mathrm{SU}_\alpha(2,1)$, see Definition (2.13). Note that if two non-degenerate Hermitian matrices are scalar multiples of each other, the corresponding unitary groups are the same. Therefore, we can conclude that case (2) parametrizes all the two dimensional unitary groups over all the quadratic extensions of the ground field.

There may be two kinds of applications of results in this paper in the future. One is of course the zeta function theory of prehomogeneous vector spaces. The above theorem says that the zeta function for these cases is the counting function for quadratic fields weighted with the Tamagawa number (without the normalization) of the group $\mathrm{SU}_\alpha(2,1)$ for case (1) and the Tamagawa number of two dimensional unitary groups. The other is an analogue of the Oppenheim conjecture for prehomogeneous vector spaces. For the zeta function theory, the reader should see [13]. For an analogue of the Oppenheim conjecture, the reader should see [12].

## §1 Rational orbits and the Galois cohomology

In this section we recall the relation between Galois cohomology and the set of rational orbits in prehomogeneous vector spaces.

We first recall the definition of Galois cohomology. Let $G$ be an algebraic group over $k$, and $k'/k$ a finite Galois extension. A 1–cocycle is a function $h = \{h_\eta\}$ from $\mathrm{Gal}(k'/k)$ to $G_{k'}$ ($h_\eta$ is the value of $h$ at $\eta \in \mathrm{Gal}(k'/k)$) satisfying the cocycle condition $h_{\eta_1 \eta_2} = h_{\eta_2} h_{\eta_1}^{\eta_2}$ for all $\eta_1, \eta_2$. If $h = \{h_\eta\}$, $i = \{i_\eta\}$ are 1–cocycles, they are equivalent if there exists $g \in G_{k'}$ such that $h_\eta = g^{-1} i_\eta g^\eta$ for all $\eta$. This defines an equivalence relation and $\mathrm{H}^1(k'/k, G)$ is the set of equivalence classes. Let $g \in G_{k'}$. We use the notation $\delta g$ for the 1–cocycle $h = \{h_\eta\}$ defined by $h_\eta = g^{-1} g^\eta$ for all $\eta \in \mathrm{Gal}(k'/k)$. The cohomology class defined by $\delta g$ does not depend on the choice of $g$ and we denote this element by 1.

We define $\mathrm{H}^1(k, G)$ to be the projective limit of $\mathrm{H}^1(k'/k, G)$ for all the finite Galois extensions $k'$. An equivalent definition is to consider continuous func-tions $\{h_\eta\}$ from $\mathrm{Gal}(\overline{k}/k)$ to $G_{\overline{k}}$ satisfying the same cocycle condition. We define $\mathrm{H}^0(k'/k, G) = \mathrm{H}^0(k, G) = G_k$. If $G$ is an abelian group, $\mathrm{H}^n(k'/k, G)$ can be defined for all $n$ and has a structure of an abelian group also.

Let

(1.1)  $$1 \to G_1 \to G_2 \to G_3 \to 1$$

be a short exact sequence of algebraic groups over $k$. This means that $G_1$ is a normal subgroup of $G_2$, the kernel of $G_2 \to G_3$ is $G_1$, and $G_{2\overline{k}} \to G_{3\overline{k}}$ is surjective. If $G_1, G_2, G_3$ are abelian and (1.1) is exact, the usual long exact sequence holds.

Let $g \in G_{3k}$. If $k'/k$ is a large enough finite Galois extension, there is an element $f \in G_{2k'}$ which maps to $g$. For a cohomology class $c$ in $\mathrm{H}^1(k, G_1)$ defined by a 1–cocycle $h = \{h_\eta\}$, we define $gc \in \mathrm{H}^1(k, G_1)$ to be the class defined by the 1–cocycle $\{fh_\eta(f^\eta)^{-1}\}$. It is known that this defines an action of $G_{3k}$ on $\mathrm{H}^1(k, G_1)$. The following lemma is an easy consequence of Proposition 38, §5.5 and Corollaire 1, §5.5 of [9].

**Lemma (1.2)** *The sequence*

$$1 \to G_{3k} \setminus \mathrm{H}^1(k, G_1) \to \mathrm{H}^1(k, G_2) \to \mathrm{H}^1(k, G_3)$$

*is exact. Moreover, if (1.1) is split, the last map is surjective.*

Note that the exactness of the sequence in (1.2) means that the inverse image of $1 \in \mathrm{H}^1(k, G_3)$ is $G_{3k} \setminus \mathrm{H}^1(k, G_1)$.

It is a familiar fact that both $\mathrm{H}^1(k, \mathrm{GL}(n))$ and $\mathrm{H}^1(k, \mathrm{SL}(n))$ are trivial. Let $k'/k$ be a finite extension. Then as remarked in [9] (see the proof of Théorème 1, §2.2) one has

$$\mathrm{H}^1(k_1, G) = \mathrm{H}^1(k, R_{k_1/k}(G))$$

for any algebraic group $G$ over $k_1$, where $R_{k_1/k}$ denotes restriction of scalars. This leads at once to the following.

**Lemma (1.3)** *Let $k_1/k$ be a finite extension, and $G = \mathrm{GL}(n)_{k_1}$ or $\mathrm{SL}(n)_{k_1}$ considered as an algebraic group over $k$. Then $\mathrm{H}^1(k, G) = \{1\}$.*

Let $(G, V)$ be an irreducible prehomogeneous vector space. For $x \in V_k^{\mathrm{ss}}$, let $g_x$ be the stabilizer and $G_x^0$ be the connected component of 1. Since $G_x$ is an algebraic group over $k$, $G_x/G_x^0$ is a finite algebraic group. Let $\mathfrak{S}_i$ be the group of permutations of $\{1, \cdots, i\}$. In the following we make the following assumption.

**Assumption (1.4)** *There exists an element $w \in V_k^{\mathrm{ss}}$ and $i$ such that*
(1) $\mathrm{Gal}(\overline{k}/k)$ *acts on* $G_w/G_w^0$ *trivially,*
(2) $G_w/G_w^0 \cong \mathfrak{S}_i$,
(3) *the sequence (1.1) splits for $G_1 = G_w^0$, $G_2 = G_w$, $G_3 = \mathfrak{S}_i$.*

All the prehomogeneous vector spaces in [11] and [5] satisfy this assumption. For prehomogeneous vector spaces in this paper, we prove this assumption in §§2,3.

For the rest of this section, we assume Assumption (1.4). With this assumption, we have the following exact sequence

(1.5)  $$1 \to \mathfrak{S}_i \setminus \mathrm{H}^1(k, G_w^0) \to \mathrm{H}^1(k, G_w) \to \mathrm{H}^1(k, \mathfrak{S}_i) \to 1.$$

We denote the surjective map $\mathrm{H}^1(k, G_w) \to \mathrm{H}^1(k, \mathfrak{S}_i)$ by $\gamma_V$. The set $\mathrm{H}^1(k, \mathfrak{S}_i)$ can be identified with the set of conjugacy classes of homomorphisms from $\mathrm{Gal}(\overline{k}/k)$ to $\mathfrak{S}_i$. The kernel of the homomorphism determines a field which belongs to $\mathfrak{Er}_i$. This defines a map $\mathrm{H}^1(k, \mathfrak{S}_i) \to \mathfrak{Er}_i$. Let $\alpha_V$ be the composition of $\gamma_V$ and the map $\mathrm{H}^1(k, \mathfrak{S}_i) \to \mathfrak{Er}_i$. If $i = 1, 2, 3$ ($i = 2$ for prehomogeneous vector spaces in this paper), $\mathrm{H}^1(k, \mathfrak{S}_i) \cong \mathfrak{Er}_i$, and the maps $\alpha_V$ and $\gamma_V$ can be identified.

If $x \in V_k^{\mathrm{ss}}$, we can choose a finite Galois extension $k'/k$ and $g \in G_{k'}$ such that $x = gw$. Then $c_x = \{g^{-1}g^\eta\}$ determines an element of $\mathrm{Ker}(\mathrm{H}^1(k, G_w) \to \mathrm{H}^1(k, G))$ (which is the set of elements which map to $1 \in \mathrm{H}^1(k, G)$). The following theorem (which does not require Assumption (1.4)) is due to Igusa [4].

**Theorem (1.6) (Igusa)** *The correspondence*

$$G_k \setminus V_x^{\mathrm{ss}} \ni x \to c_x \in \mathrm{Ker}(\mathrm{H}^1(k, G_w) \to \mathrm{H}^1(k, G))$$

*is bijective.*

By the above theorem, we can consider $G_k \setminus V_k^{\mathrm{ss}}$ as a subset of $\mathrm{H}^1(k, G_w)$. We denote the restriction of $\alpha_V, \gamma_V$ to $G_k \setminus V_k^{\mathrm{ss}}$ also by $\alpha_V, \gamma_V$.

Let $x \in V_k^{\mathrm{ss}}$. We choose an element $g_x \in G_{\overline{k}}$ so that $x = g_x w$. Then for each element $c \in \mathrm{H}^1(k, G_w)$ defined by a 1–cocycle $\{h_\eta\}$, we can associate an element $c^{g_x} \in \mathrm{H}^1(k, G_x)$ defined by a 1–cocycle $\{g_x h_\eta (g_x^\eta)^{-1}\}$. Note that since $(g_x^\eta)^{-1} g_x \in G_w$ (but may not be in $G_w^0$),

$$g_x h_\eta (g_x^\eta)^{-1} = g_x (h_\eta (g_x^\eta)^{-1} g_x) g_x^{-1} \in G_{x\overline{k}}.$$

It is easy to see that the map $\mathrm{H}^1(k, G_w) \ni c \to c^{g_x} \in \mathrm{H}^1(k, G_x)$ is well defined and does not depend on the choice of $g_x$. Also a similar construction using $g_x^{-1}$ defines a map from $\mathrm{H}^1(k, G_x)$ to $\mathrm{H}^1(k, G_w)$. So $c \to c^{g_x}$ is a bijective map.

Let $x \in V_k^{\mathrm{ss}}$. By the exact sequence

(1.7) $$1 \to G_x^0 \to G_x \to G_x/G_x^0 \to 1$$

and Lemma (1.2), $(G_x/G_x^0)_k \setminus \mathrm{H}^1(k, G_x^0)$ can be considered as a subset of $\mathrm{H}^1(k, G_x)$, which corresponds bijectively with $\mathrm{H}^1(k, G_w)$. It is easy to see that $c \in \mathrm{H}^1(k, G_w)$ maps to the trivial element of $\mathrm{H}^1(k, G)$ if and only if $c^{g_x} \in \mathrm{H}^1(k, G_x)$ maps to the trivial element of $\mathrm{H}^1(k, G)$ by the imbedding $G_x \subset G$. So cohomology classes in $(G_x/G_x^0)_k \setminus \mathrm{H}^1(k, G_x^0)$ which correspond to rational orbits are those in

$$\mathrm{Ker}((G_x/G_x^0)_k \setminus \mathrm{H}^1(k, G_x^0) \to \mathrm{H}^1(k, G)).$$

The proof of Lemma (1.12) [5] works for the following Lemma also and we do not repeat the proof here. Note that $\overline{k} = k^{\mathrm{sep}}$ because we assumed ch $k = 0$. We used the notation $\gamma_V$ here instead of $\alpha_V$ in [5] for the map $G_k \setminus V_k^{\mathrm{ss}} \to \mathrm{H}^1(k, \mathfrak{S}_i)$ because there is a slight difference between $\mathrm{H}^1(k, \mathfrak{S}_i)$ and $\mathfrak{Er}_i$ if $i \geq 4$, but that was not the case in [5] because $i = 2, 3$ in [5].

**Lemma (1.8)** (1) *If* $x \in V_k^{\mathrm{ss}}$,

$$\gamma_V^{-1}(\gamma_V(G_k x)) \cong \mathrm{Ker}((G_x/G_x^0)_k \setminus \mathrm{H}^1(k, G_x^0) \to \mathrm{H}^1(k, G)).$$

(2) *By this identification, the cohomology class* $\{g^{-1}g^\eta\} \in \mathrm{H}^1(k, G_x^0)$ $(g \in G_{\bar{k}})$ *corresponds to the orbit* $G_k g x$.

## §2 The orbit space $G_k \setminus V_k^{\mathrm{ss}}$ (1)

Let $G_1, W, c$, etc., be as in the introduction. In this section, we consider the prehomogeneous vector space $G = G_1$, $V = W$. We prove that $\alpha_V$ for this case is a map to $\mathfrak{Er}_2$ and construct a subset of $G_k \setminus V_k^{\mathrm{ss}}$ which maps bijectively to $\mathfrak{Er}_2$. Also for each point $x$ in this set, we determine the stabilizer $G_x$.

Let $x = {}^t(x_1, \cdots, x_7)$ be the coordinate of $x \in V$ with respect to $\{e_1, \cdots, e_7\}$. It is known (see [8, p. 135]) that $(G, V)$ is a prehomogeneous vector space and

$$(2.1) \qquad \Delta(x) = x_1^2 - 4(x_2 x_5 + x_3 x_6 + x_4 x_7)$$

is a relative invariant polynomial. It is proved in [10] that $V^*$ can be considered as the imaginary part of the split octonion algebra (which is also called the split Cayley algebra). By $\Delta$, we can identify $V$ with $V^*$ and we can regard $\Delta$ as a function on $V^*$. The resulting function on $V^*$ turns out to be a scalar multiplication of the restriction of the norm form of the octonion algebra to the imaginary part. Moreover, if $\{f_1, \cdots, f_7\}$ is the dual basis of $\{e_1, \cdots, e_7\}$ and $y = y_1 f_1 + \cdots + y_7 f_7$, $\Delta(y) = y_1^2 - (y_2 y_5 + y_3 y_6 + y_4 y_7)$.

If $k$ is algebraically closed, $c$ is set-theoretically surjective and the kernel of $c$ is a semi-simple group of type $\mathrm{G}_2$ which is connected modulo scalar matrices (for the second part, see [8, p. 85]). So $G$ is a connected reductive group over any ground field $k$ of characteristic zero. If $\Delta(gx) = \chi(g)\Delta(x)$, $c(tI_7) = t^3$, $\chi(tI_7) = t^2$. Since $c, \chi$ are rational characters, so is $\chi' = c\chi^{-1}$ and $\chi'(tI_7) = t$. Since $G$ is connected and $\dim G/[G, G] = 1$, the group of rational characters is generated by one element. Therefore, $c = \chi'^3$, $\chi = \chi'^2$.

Let $w = {}^t(1 \ 0 \ \cdots \ 0)$. Then $w \in V_k^{\mathrm{ss}}$ because $\Delta(w) = 1 \neq 0$. It is proved in [8, p. 134] that the Lie algebra of $G_w^0$ is

$$\left\{ \begin{pmatrix} 1 & 0 & 0 \\ 0 & A & 0 \\ 0 & 0 & -{}^t A \end{pmatrix} \;\middle|\; A \in \mathrm{M}(3,3),\ \mathrm{tr}(A) = 0 \right\}.$$

Consider elements of the form

$$(2.2) \qquad g(a, b, c, A) = g(a, b_1, b_2, c_1, c_2, A_1, A_2, A_3, A_4) = \begin{pmatrix} a & b_1 & b_2 \\ c_1 & A_1 & A_2 \\ c_2 & A_3 & A_4 \end{pmatrix}$$

for $a \in \mathrm{M}(1,1)$, $b_1, b_2 \in \mathrm{M}(1,3)$, $c_1, c_2 \in \mathrm{M}(3,1)$, $A_1, \cdots, A_4 \in \mathrm{M}(3,3)$. Let

$$F = \{g(a, b, c, A) \in G \mid a = 1,\ b = 0,\ c = 0,\ A_2 = A_3 = 0\}.$$

Then $F$ is clearly a closed subset of $G$. For $A \in \mathrm{SL}(3)$, let

$$(2.3) \qquad d(A) = g(1, 0, 0, 0, 0, A, 0, {}^t A^{-1}, 0).$$

It is easy to see that if $g(a, b, c, A) \in F$, $A_1{}^t A_4$ is a scalar matrix, say $dI_3$, and $\det A_1 = \det A_4 = d$. This implies that $\det A_1 = \det A_4 = d = 1$. Therefore, we get the following lemma.

**Lemma (2.4)** $F = \{d(A) \mid A \in \mathrm{SL}(3)\}$.

The above lemma implies that $\mathrm{SL}(3)$ is a closed subgroup of $G$ by the imbedding $A \to d(A)$. Clearly, elements of the form $d(A)$ fix the point $w$. Since $\mathrm{SL}(3)$ is connected, $\mathrm{SL}(3) \subset G_w^0$ is a closed subgroup. Since they have the same Lie algebra, $G_{w\overline{k}}^0 = \mathrm{SL}(3)_{\overline{k}}$ (see [2]). Therefore, by Lemma 2 [6, p. 8], $\mathrm{SL}(3) = G_w^0$ as algebraic groups over $k$. In particular, $\mathrm{SL}(3)$ is a normal subgroup of $G_w$.

Let

$$(2.5) \qquad \Lambda = \begin{pmatrix} 1 & & \\ & 1 & \\ & & -1 \end{pmatrix}, \ \tau = \begin{pmatrix} 1 & & \\ & & \Lambda \\ & \Lambda & \end{pmatrix}.$$

**Proposition (2.6)** *The subgroup $G_w$ is generated by $G_w^0 \cong \mathrm{SL}(3)$ and $\tau$.*

*Proof.* We may assume that $k$ is algebraically closed by Lemma 2 [6, p. 8] again. Note that $\tau^2 = 1$ and $\tau$ induces an outer automorphism of $\mathrm{SL}(3)$ by conjugation. This implies that conjugations by elements of the group generated by $G_w^0$ and $\tau$ exhaust all the automorphisms of $\mathrm{SL}(3)$. Let $g \in G_w$. Then by the above observation, multiplying an element of $\mathrm{SL}(3)$ and $\tau$ if necessary, we may assume that $g$ commutes with elements of $\mathrm{SL}(3)$.

We define

$$(2.7) \qquad U_0 = \langle e_1 \rangle, \ U_1 = \langle e_2, e_3, e_4 \rangle, \ U_2 = \langle e_5, e_6, e_7 \rangle$$

($\langle \ \rangle$ means the span). Then $U_0$ is the trivial representation of $\mathrm{SL}(3)$, $U_1$ is the standard representation of $\mathrm{SL}(3)$, and $U_2$ is the dual of $U_1$. As representations of $\mathrm{SL}(3)$, $V = U_0 \oplus U_1 \oplus U_2$, and $U_0, U_1, U_2$ are inequivalent representations of $\mathrm{SL}(3)$. Since $g$ commutes with elements of $G_w^0$, $g(U_i) \subset U_i$ for $i = 0, 1, 2$ and by Schur's lemma, $g$ must be of the form $g(a, 0, 0, 0, 0, bI_3, 0, cI_3, 0)$. It is easy to verify that this element fixes $w$ if and only if $a = b = c = 1$. This proves the proposition. $\square$

By the above proposition, we have a split exact sequence

$$(2.8) \qquad 1 \to \mathrm{SL}(3) \to G_w \to \mathbb{Z}/2\mathbb{Z} \to 1,$$

where $\tau$ maps to the non-trivial element of $\mathbb{Z}/2\mathbb{Z}$.

Because of the split exact sequence (2.8), $\mathrm{H}^1(k, \mathbb{Z}/2\mathbb{Z}) \cong \mathfrak{Ex}_2$ can be considered as a subset of $\mathrm{H}^1(k, G_w)$.

**Lemma (2.9)** *The map $\mathfrak{Ex}_2 \to \mathrm{H}^1(k, G)$ is trivial.*

*Proof.* The field $k$ as the degree one extension of $k$ corresponds to the trivial class of $\mathrm{H}^1(k, \mathbb{Z}/2\mathbb{Z})$ and the lemma is trivial for this case. Let $k(\alpha)/k$ be a quadratic extension such that $\alpha \notin k$ and $\alpha^2 \in k$. Let $\sigma$ be the non-trivial element of $\mathrm{Gal}(k(\alpha)/k)$.

Let
$$h_0 = \begin{pmatrix} 0 & -1 & 0 & 0 & 1 & 0 & 0 \\ \frac{1}{2} & \frac{1}{2} & 0 & 0 & \frac{1}{2} & 0 & 0 \\ 0 & 0 & \frac{1}{2} & -\frac{1}{2} & 0 & -\frac{1}{2} & -\frac{1}{2} \\ 0 & 0 & \frac{1}{2} & \frac{1}{2} & 0 & \frac{1}{2} & -\frac{1}{2} \\ -\frac{1}{2} & \frac{1}{2} & 0 & 0 & \frac{1}{2} & 0 & 0 \\ 0 & 0 & -\frac{1}{2} & -\frac{1}{2} & 0 & \frac{1}{2} & -\frac{1}{2} \\ 0 & 0 & \frac{1}{2} & -\frac{1}{2} & 0 & \frac{1}{2} & \frac{1}{2} \end{pmatrix}.$$

Then easy computations show that $h_0 \in G$ (in fact $h_0 \in G \cap \mathrm{SL}(7)$).

We define

(2.10) $$g_\alpha = \alpha I_7 \begin{pmatrix} 1 & & & & & & \\ & \alpha & & & & & \\ & & 1 & & & & \\ & & & \alpha^{-1} & & & \\ & & & & \alpha^{-1} & & \\ & & & & & 1 & \\ & & & & & & \alpha \end{pmatrix} h_0.$$

Then
$$g_\alpha = \begin{pmatrix} 0 & -\alpha & 0 & 0 & \alpha & 0 & 0 \\ \frac{\alpha^2}{2} & \frac{\alpha^2}{2} & 0 & 0 & \frac{\alpha^2}{2} & 0 & 0 \\ 0 & 0 & \frac{\alpha}{2} & -\frac{\alpha}{2} & 0 & -\frac{\alpha}{2} & -\frac{\alpha}{2} \\ 0 & 0 & \frac{1}{2} & \frac{1}{2} & 0 & \frac{1}{2} & -\frac{1}{2} \\ -\frac{1}{2} & \frac{1}{2} & 0 & 0 & \frac{1}{2} & 0 & 0 \\ 0 & 0 & -\frac{\alpha}{2} & -\frac{\alpha}{2} & 0 & \frac{\alpha}{2} & -\frac{\alpha}{2} \\ 0 & 0 & \frac{\alpha^2}{2} & -\frac{\alpha^2}{2} & 0 & \frac{\alpha^2}{2} & \frac{\alpha^2}{2} \end{pmatrix}.$$

This implies that $g^\sigma = g_\alpha \tau$. Therefore, the cohomology class in $\mathrm{H}^1(k, \mathbb{Z}/2\mathbb{Z})$ which corresponds to the field $k(a)$ maps to the trivial class $\{g_\alpha^{-1} g_\alpha^\eta\}_{\eta \in \mathrm{Gal}(\overline{k}/k)}$ in $\mathrm{H}^1(k, G)$. □

**Definition (2.11)** $w_\alpha = g_\alpha w = {}^t\begin{pmatrix} 0 & \frac{\alpha^2}{2} & 0 & 0 & -\frac{1}{2} & 0 & 0 \end{pmatrix}$.

By the correspondence in Theorem (1.7), $w_\alpha \in V_k^{\mathrm{ss}}$ and it corresponds to the field $k(\alpha)$. Since $\Delta(w_\alpha) = \alpha^2$, we get the following proposition.

**Proposition (2.12)** $k(\alpha) = k(\Delta(w_\alpha)^{\frac{1}{2}})$.

Let $R$ be any $k$-algebra. We define $R(\alpha) = R \otimes k(\alpha)$. The non-trivial element $\sigma \in \mathrm{Gal}(k(\alpha)/k)$ acts on $R(\alpha)$ by $(r \otimes x)^\sigma = r \otimes (x^\sigma)$. For $A \in \mathrm{SL}(3)_{R(\alpha)}$, we define $A^\sigma, {}^*A = {}^tA^\sigma$ entry-wise. Let

(2.13) $$\mathrm{SU}_\alpha(2,1)_R = \{A \in \mathrm{SL}(3)_{R(\alpha)} \mid A\Lambda {}^*A = \Lambda\}.$$

We define $\mathrm{SU}_\alpha(2,1)$ to be the group over $k$ such that the set of $R$–valued points is $\mathrm{SU}_\alpha(2,1)_R$ for any $k$–algebra $R$.

**Proposition (2.14)** $G^0_{w_\alpha} = \{g_\alpha d(A) g_\alpha^{-1} \mid A \in \mathrm{SU}_\alpha(2,1)\} \cong \mathrm{SU}_\alpha(2,1)$.

*Proof.* Let $R$ be any $k$–algebra. Consider $R(\alpha)$ as above with the action of $\sigma$. Regarding $G_{w_\alpha} \subset \mathrm{GL}(7)$, we consider the entry-wise action of $\sigma$ for elements in $G^0_{w_\alpha R(\alpha)}$ also. Then
$$G^0_{w_\alpha R} = \{g \in G^0_{w_\alpha R(\alpha)} \mid g^\sigma = g\}.$$

Over $R(\alpha)$, we can express elements of $G^0_{w_\alpha R(\alpha)}$ as $g_\alpha d(A) g_\alpha^{-1}$ where $A \in \mathrm{SL}(3)_{R(\alpha)}$. This element is in $G^0_{w_\alpha R}$ if and only if

$$\begin{aligned} g_\alpha d(A) g_\alpha^{-1} &= g_\alpha^\sigma d(A^\sigma)(g_\alpha^\sigma)^{-1} \\ &= g_\alpha \tau d(A^\sigma) \tau g_\alpha^{-1} \\ &= g_\alpha d(\Lambda^* A^{-1} \Lambda) g_\alpha^{-1}. \end{aligned}$$

This condition is satisfied if and only if $A \in \mathrm{SU}_\alpha(2,1)_R$. Therefore, the proposition follows from Theorem [7, p. 17]. $\square$

### §3 The orbit space $G_k \backslash V_k^{\mathrm{ss}}$ (2)

In this section, we consider case (2) in the introduction. Let $G_1, W, c$, etc. be as in the introduction, and $\chi$ as in §2. We define the prehomogeneous vector space (2) more precisely first.

Let $G_2 = \mathrm{GL}(2)$, and $G = G_1 \times G_2$. Let $V = W \otimes k^2$ and we consider this space as the set of linear forms $x = M(v)$ in two variables $v = (v_1, v_2)$ with entries in $W$. Let $W_1 = v_1 W$ be the subspace of forms $x = M(v)$ such that $M(0, v_2) = 0$. We define $W_2 = v_2 W$ similarly. Then $W = W_1 \oplus W_2$.

We define an action of $g = (g_1, g_2) \in G$ on $x = M(v)$ by $gx = g_1 M(vg_2)$. This defines a representation of $G$ on $V$. Let $\{e_1, \cdots, e_7\}$ be the standard basis of $W$. Considering $V = v_1 W \oplus v_2 W$, let $e_{i,j} = v_i e_j$ for $i = 1, 2, j = 1, \cdots, 7$. Then $\{e_{i,j} \mid i = 1, 2, j = 1, \cdots, 7\}$ is a basis of $V$. Let $x = (x_{i,j})$ be the coordinate of $x \in V$ with respect to this basis. We may also use the notation $x = (x_1, x_2)$ where $x_i = {}^t(\,x_{i,1} \;\cdots\; x_{i,7}\,)$ for $i = 1, 2$.

Let $\widetilde{T} = \mathrm{Ker}(G \to \mathrm{GL}(V))$. Then
$$\widetilde{T} = \{(t^{-1} I_7, t I_2) \mid t \in \mathrm{GL}(1)\} \cong \mathrm{GL}(1),$$

and $\widetilde{T}$ is contained in the center of $G$.

Consider the element

(3.1) $$w = v_1 \begin{pmatrix} 0 \\ 1 \\ 0 \\ 0 \\ 0 \\ 0 \\ 0 \end{pmatrix} + v_2 \begin{pmatrix} 0 \\ 0 \\ 0 \\ 0 \\ 1 \\ 0 \\ 0 \end{pmatrix}.$$

Let $\Delta$ be as in (2.1). For $x = (x_1, x_2) \in V$, we define

(3.2) $$F_x(v) = \Delta(v_1 x_1 + v_2 x_2).$$

This is a binary quadratic form. It is easy to see that if $g = (g_1, g_2) \in G$,

$$F_{gx}(v) = \chi(g_1) F_x(v g_2).$$

Therefore, the discriminant of $F_x(v)$ is a relative invariant polynomial. Note that $F_w(v) = -4 v_1 v_2$. So $w \in V_k^{ss}$.

**Definition (3.3)** *For $A \in \mathrm{GL}(2)$, we define $d(A) = (d_1(A), d_2(A))$ where*

$$d_1(A) = \begin{pmatrix} (\det A)^{-1} & & & \\ & A & & \\ & & \det A & \\ & & & {}^t A^{-1} \end{pmatrix}, \quad d_2(A) = \begin{pmatrix} \det A & \\ & (\det A)^{-1} \end{pmatrix}.$$

By the imbedding $A \to \begin{pmatrix} (\det A)^{-1} & 0 \\ 0 & A \end{pmatrix}$, $\mathrm{GL}(2)$ becomes a closed subgroup of $\mathrm{SL}(3) \subset G_1$. So $\mathrm{GL}(2)$ is a closed subgroup of $G$ by the imbedding $A \to d(A)$. It is easy to see that elements of the form $d(A)$ fix $w$. Since $\mathrm{GL}(2)$ is connected, $\mathrm{GL}(2)$ is a closed subgroup of $G_w^0$. It is known (see [8, p. 136]) that $\mathrm{GL}(2)$ and $G_w^0 / \widetilde{T}$ have the same Lie algebra. It is easy to see that $\mathrm{GL}(2) \cap \widetilde{T} = \{1\}$. Therefore, $G_w^0 \cong \mathrm{GL}(2) \times \widetilde{T}$. In particular $\mathrm{GL}(2) \times \widetilde{T}$ is a normal subgroup of $G_w$.

We define $\tau = (\tau_1, \tau_2)$ where

$$(3.4) \qquad \Lambda = \begin{pmatrix} 1 & \\ & -1 \end{pmatrix}, \quad \tau_1 = \begin{pmatrix} & & 1 & \\ & 1 & & \\ & & & \Lambda \\ 1 & & & \\ & & \Lambda & \end{pmatrix}, \quad \tau_2 = \begin{pmatrix} & 1 \\ 1 & \end{pmatrix}.$$

Note that we are using different definitions for $d(A), \Lambda, \tau$ in this section ($\tau_1$ is $\tau$ in §2.)

**Proposition (3.5)** *The subgroup $G_w$ is generated by $G_w^0 \cong \mathrm{GL}(2) \times \widetilde{T}$ and $\tau$.*

*Proof.* As in §2, we may assume that $k$ is algebraically closed. Suppose $U$ is an irreducible representation of $\mathrm{GL}(2)$. Then there exist two integers $k, l$ such that the highest weight of elements of the form $\begin{pmatrix} t & 0 \\ 0 & t^{-1} \end{pmatrix}$ is $t^k$ and $t_2 I_2$ acts by $t_2^l$. Note that $k \geq 0$. Since $U$ is determined by $k, l$, we say $U$ is of type $(k, l)$.

Let

$$(3.6) \qquad \begin{aligned} U_{i,0} &= \langle e_{i,1} \rangle, \ U_{i,1} = \langle e_{i,2} \rangle, \ U_{i,2} = \langle e_{i,3}, e_{i,4} \rangle, \\ U_{i,3} &= \langle e_{i,5} \rangle, \ U_{i,4} = \langle e_{i,6}, e_{i,7} \rangle \end{aligned}$$

for $i = 1, 2$.

Suppose $g = (g_1, g_2) \in G_w$. Then $F_w(v g_2)$ is a scalar multiple of $F_w(v)$. Multiplying $\tau$ if necessary, we may assume that $(1 \ \ 0) g_2$, $(0 \ \ 1) g_2$ are scalar multiples of $(1\ 0)$, $(0\ 1)$ respectively. This implies $g_2$ is a diagonal element and commutes with $d_2(A)$ for any $A \in \mathrm{GL}(2)$.

Since $g$ normalizes $\mathrm{GL}(2) \times \widetilde{T}$, the conjugation by $g$ maps the center of $\mathrm{GL}(2) \times \widetilde{T}$ to itself. Since $\widetilde{T}$ is contained in the center of $G$, the conjugation by $g_1$ induces either the identity map or the inverse on the center of $\mathrm{GL}(2) \subset G_1$ (imbedded by $d_1$) modulo $\widetilde{T}$. Since $g_2$ commutes with $d_1(A)$ for any $A \in \mathrm{GL}(2)$, $g_1$ must induce the identity map on the center of $\mathrm{GL}(2) \subset G_1$ modulo $\widetilde{T}$. Therefore, $g$ commutes with the center of $\mathrm{GL}(2) \subset G$ modulo $\widetilde{T}$.

Since $g_2$ is diagonal, $g(W_1) = W_1$, $g(W_2) = W_2$. Note that $\widetilde{T}$ acts trvially on $V$. So by the above consideration, $g$ does not change the weights with respect to $tI_2$ of elements in $W_1$. The types of $U_{1,j}$ for $j = 0, \cdots, 4$ are

$$(0,2), (0,0), (1,3), (0,4), (1,1).$$

Therefore, $g_1$ must be of the form

$$g_1 = \begin{pmatrix} a & & & & \\ & b & & & \\ & & A & & \\ & & & b^{-1} & \\ & & & & {}^tA^{-1} \end{pmatrix},$$

where $a, b \in \mathrm{GL}(1)$, $A \in \mathrm{GL}(2)$.

Multiplying $d(bA^{-1})$ and an element of $\widetilde{T}$, we may assume that $b = 1$, $A = I_2$. It is easy to see that if

$$g = \left( \begin{pmatrix} a & 0 \\ 0 & I_6 \end{pmatrix}, \begin{pmatrix} b & 0 \\ 0 & c \end{pmatrix} \right) \in G$$

fixes $w$, $a = b = c = 1$. This proves the proposition. $\square$

By the above proposition, we have a split exact sequence

(3.7) $$1 \to \mathrm{GL}(2) \times \mathrm{GL}(1) \to G_w \to \mathbb{Z}/2\mathbb{Z} \to 1$$

where $\tau$ maps to the non-trivial element of $\mathbb{Z}/2\mathbb{Z}$. As in §2,

$$G_k \backslash V_k^{\mathrm{ss}} \cong \mathrm{Ker}(\mathrm{H}^1(k, G_w) \to \mathrm{H}^1(k, G))$$

and we have a map $\mathrm{H}^1(k, G_w) \supset \mathfrak{Er}_2 \to \mathrm{H}^1(k, G)$.

Let $k(\alpha)/k$ be a quadratic extension such that $\alpha \notin k$ and $\alpha^2 \in k$. Let $g_{\alpha,1}$ be the matrix defined in (2.10) as $g_\alpha$. Let $g_{\alpha,2} = \begin{pmatrix} 1 & 1 \\ \alpha & -\alpha \end{pmatrix}$, and $g_\alpha = (g_{\alpha,1}, g_{\alpha,2}) \in G_{k(\alpha)}$.

**Definition (3.8)** $w_\alpha = g_\alpha w$.

Let $\sigma \in \mathrm{Gal}(k(\alpha)/k)$ the non-trivial element as before. Since $g_\alpha$ satisfies the condition $g_\alpha^\sigma = g_\alpha \tau$, we get the following proposition.

**Proposition (3.9)** *The map $\mathfrak{Ex}_2 \to \mathrm{H}^1(k, G)$ is trivial, $w_\alpha \in V_k^{\mathrm{ss}}$, and $w_\alpha$ corresponds to the field $k(\alpha)$.*

Since the roots of the polynomial $F_{w_\alpha}(v)$ are $\pm \alpha$, we get the following corollary.

**Corollary (3.10)** *The field $k(\alpha)$ is generated by roots of the polynomial $F_{w_\alpha}(v)$.*

For $A \in \mathrm{GL}(2)_{k(\alpha)}$, we define $^*A$ as in §2.

**Definition (3.11)** $\mathrm{U}_\alpha(1,1) = \{A \in \mathrm{GL}(2)_{k(\alpha)} \mid A\Lambda^* A = \Lambda\}$.

We consider $\mathrm{U}_\alpha(1,1)$ as an algebraic group over $k$. The proof of the following proposition is similar to that of Proposition (2.15) and is left to the reader.

**Proposition (3.12)** $G^0_{w_\alpha} = \left\{g_\alpha d(A) g_\alpha^{-1} \,\middle|\, A \in \mathrm{U}_\alpha(1,1)\right\} \times \widetilde{T} \cong \mathrm{U}_\alpha(1,1) \times \widetilde{T}$.

### §4 The fiber structure of $\alpha_V$

We consider the fiber structure of $\alpha_V$ for prehomogeneous vector spaces (1), (2) in this section. In both cases, $\alpha_V$ is a map to $\mathfrak{Ex}_2$, there is a point $w \in V_k^{\mathrm{ss}}$ which corresponds to the trivial extension $k/k$, and $G_w^0 \cong \mathrm{SL}(3)$ or $\mathrm{GL}(2) \times \mathrm{GL}(1)$. Since $\mathrm{H}^1(k, \mathrm{SL}(3)) = \mathrm{H}^1(k, \mathrm{GL}(2) \times \mathrm{GL}(1)) = \{1\}$, $\alpha_V^{-1}(k) = G_k w$ in both cases by Lemma (1.8).

Let $k(\alpha)/k$ be a quadratic extension such that $\alpha^2 \in k$ and $\sigma \in \mathrm{Gal}(k(\alpha)/k)$ is the non-trivial element throughout this section. We proved in §§2,3 that $\alpha_V$ is surjective and chose a point $w_\alpha = g_\alpha w$ which corresponds to the field $k(\alpha)$ in both cases. To determine the fiber over $k(\alpha)/k$, we consider the map $\mathrm{H}^1(k, G^0_{w_\alpha}) \to \mathrm{H}^1(k, G)$.

We first consider case (1). So we use the definitions of $G, V, w, d(A), \Lambda, \tau, g_\alpha, w_\alpha$ in §2 in (4.1)–(4.5).

Let $F_1 = \mathrm{SU}_\alpha(2,1)$ and $F_2$ be the restriction of $\mathrm{SL}(3)$ from $k(\alpha)$ to $k$. For $g \in \mathrm{SL}(3)_{k(\alpha)}$, we denote the entry-wise action of $\sigma$ by $g^\sigma$. We define $^*g = {}^t g^\sigma$. For $g \in F_{ik(\alpha)}$ ($i = 1, 2$), let $\phi_{F_i}(\sigma)(g)$ be the action of $\sigma$ induced from the $k$–structure of $F_i$.

We use the notation $i_k, i_{k(\alpha)}$ for the natural imbeddings $F_{1k} \to F_{2k}, F_{1k(\alpha)} \to F_{2k(\alpha)}$. We denote the natural inclusion map $F_{ik} \to F_{ik(\alpha)}$ by $j_{ik(\alpha)}$ for $i = 1, 2$. Obviously, the following diagram is commutative.

$$\begin{array}{ccc} F_{1k} & \overset{j_{1k(\alpha)}}{\to} & F_{1k(\alpha)} \\ i_k \downarrow & & \downarrow i_{k(\alpha)} \\ F_{2k} & \overset{j_{2k(\alpha)}}{\to} & F_{2k(\alpha)} \end{array}$$

We can identify $F_{1k(\alpha)} \cong \mathrm{SL}(3)_{k(\alpha)}$, $F_{2k(\alpha)} \cong \mathrm{SL}(3)_{k(\alpha)} \times \mathrm{SL}(3)_{k(\alpha)}$. With this identification,

$$\begin{aligned}(4.1) \quad & i_{k(\alpha)}(g) = (g, \Lambda {}^t g^{-1} \Lambda), \; j_{1k(\alpha)}(g) = g, \; j_{2k(\alpha)}(g) = (g, g^\sigma), i_k(g) = g, \\ & \phi_{F_1}(\sigma)(g) = \Lambda^* g^{-1} \Lambda, \; \phi_{F_2}(\sigma)(g_1, g_2) = (g_2^\sigma, g_1^\sigma). \end{aligned}$$

For $i = 1, 2$, we can consider $g \to g^\sigma$ as a $k$–automorphism of the $k$–group $F_i$ and denote it by $\psi_{F_i}(\sigma)$. Note that $\psi_{F_2}(\sigma)(g_1, g_2) = (g_2, g_1)$ for $(g_1, g_2) \in F_{2k(\alpha)}$.

Let $U$ be the $k$–vector space of $3 \times 3$ Hermitian matrices (with respect to the extension $k(\alpha)/k$). Then $F_2$ acts on $U$ by $H \to gH {}^t\psi_{F_2}(\sigma)(g)$. We can identify $U_{k(\alpha)} \cong \mathrm{M}(3,3)_{k(\alpha)}$ and the action of $(g_1, g_2) \in F_{2k(\alpha)}$ is $H \to g_1 H {}^t g_2$. The determinant $\det H$ is a function on $U$ defined over $k$. We define $X = \{H \in U \mid \det H = -1\}$.

Since $\operatorname{Gal}(k(\alpha)/k) \cong \mathbb{Z}/2\mathbb{Z}$, any element of $\mathrm{H}^1(k(\alpha)/k, F_1)$ is determined by a single element $h \in F_{1k(\alpha)} \cong \operatorname{SL}(3)_{k(\alpha)}$ satisfying the cocycle condition $h\phi_{F_1}(\sigma)(h) = h\Lambda^* h^{-1}\Lambda = 1$. So we express such a cohomology class by $\{h\}$. Note that the cocycle condition for $\{h\}$ is equivalent to $h\Lambda$ or $h^{-1}\Lambda$ being Hermitian.

**Proposition (4.2)** (1) $G_k \backslash X_k \cong \mathrm{H}^1(k, F_1)$.

(2) *Any element of* $\mathrm{H}^1(k, F_1)$ *is represented by an element of* $\mathrm{H}^1(k(\alpha)/k, F_1)$.

(3) *If* $\{h\} \in \mathrm{H}^1(k(\alpha)/k, F_1)$, *the corresponding Hermitian form is* $h^{-1}\Lambda$.

*Proof.* For $H \in X_k$, we can find an element $g \in F_{2\overline{k}}$ which maps $\Lambda$ to $H$. Then $c_H = \{g^{-1}g^\eta\}_{\eta \in \operatorname{Gal}(\overline{k}/k)}$ defines a cohomology class in $\mathrm{H}^1(k, F_1)$. Note that $\mathrm{H}^1(k, F_2) = \{1\}$. So by the argument of Igusa [4], the map $G_k \backslash X_k \ni H \to c_H \in \mathrm{H}^1(k, F_1)$ is bijective. This proves (1).

For $H \in X_k$, put $g_1 = H\Lambda^{-1}$. Then $\det g_1 = 1$. So if we choose $g = (g_1, 1)$, $g \in F_{2k(\alpha)}$ and it maps $\Lambda$ to $H$. This means that $c_H$ is represented by the element $\{g^{-1}g^\sigma\} \in \mathrm{H}^1(k(\alpha)/k, F_1)$. This proves (2). Also if we express $c_H = \{h\}$,

$$\begin{aligned} i_{k(\alpha)}(h) &= (g_1, 1)^{-1}\phi_{F_2}(\sigma)(g_1, 1) \\ &= (g_1^{-1}, g_1^\sigma) \\ &= (g_1, \Lambda g_1 \Lambda) \text{ since } H \text{ is Hermitian} \\ &= i_{k(\alpha)}(g_1^{-1}). \end{aligned}$$

Therefore, $H = g_1\Lambda = h^{-1}\Lambda$. This proves (3). $\square$

**Lemma (4.3)** *If* $x \in V_k^{\mathrm{ss}}$, *there exists* $g \in G_k \cap \operatorname{SL}(7)_k$ *such that* $gx$ *is of the form* $(0 \ 1 \ 0 \ 0 \ x_5 \ 0 \ 0)$.

*Proof.* Let

$$u_1(a,b,c) = \exp\left(\begin{pmatrix} 0 & 0 & 0 & 0 & 2a & 2b & 2c \\ a & 0 & 0 & 0 & 0 & 0 & 0 \\ b & 0 & 0 & 0 & 0 & 0 & 0 \\ c & 0 & 0 & 0 & 0 & 0 & 0 \\ 0 & 0 & -c & b & 0 & 0 & 0 \\ 0 & c & 0 & -a & 0 & 0 & 0 \\ 0 & -b & a & 0 & 0 & 0 & 0 \end{pmatrix}\right)$$

$$= \begin{pmatrix} 1 & 0 & 0 & 0 & 2a & 2b & 2c \\ a & 1 & 0 & 0 & a^2 & ab & ac \\ b & 0 & 1 & 0 & ab & b^2 & bc \\ c & 0 & 0 & 1 & ac & bc & c^2 \\ 0 & 0 & -c & b & 1 & 0 & 0 \\ 0 & c & 0 & -a & 0 & 1 & 0 \\ 0 & -b & a & 0 & 0 & 0 & 1 \end{pmatrix},$$

$$u_2(d,e,f) = \exp\left(\begin{pmatrix} 0 & 2d & 2e & 2f & 0 & 0 & 0 \\ 0 & 0 & 0 & 0 & 0 & f & -e \\ 0 & 0 & 0 & 0 & -f & 0 & d \\ 0 & 0 & 0 & 0 & e & -d & 0 \\ d & 0 & 0 & 0 & 0 & 0 & 0 \\ e & 0 & 0 & 0 & 0 & 0 & 0 \\ f & 0 & 0 & 0 & 0 & 0 & 0 \end{pmatrix}\right)$$

$$= \begin{pmatrix} 1 & 2d & 2e & 2f & 0 & 0 & 0 \\ 0 & 1 & 0 & 0 & 0 & f & -e \\ 0 & 0 & 1 & 0 & -f & 0 & d \\ 0 & 0 & 0 & 1 & e & -d & 0 \\ d & d^2 & de & df & 1 & 0 & 0 \\ e & de & e^2 & ef & 0 & 1 & 0 \\ f & df & ef & f^2 & 0 & 0 & 1 \end{pmatrix}.$$

Then $u_1(a,b,c), u_2(d,e,f) \in G_k \cap \mathrm{SL}(7)_k$. Let $x = {}^t(x_1 \cdots x_7)$ and $\Delta(x) \neq 0$. Suppose $(x_2, \cdots, x_7) \neq 0$. Then $(x_2, x_3, x_4) \neq 0$ or $(x_5, x_6, x_7) \neq 0$. Multiplying $\tau$ if necessary, we may assume that $(x_2, x_3, x_4) \neq 0$. Multiplying an element of $\mathrm{SL}(3)$, we may assume that $(x_2, x_3, x_4) = (1, 0, 0)$. Then $u_2(-\frac{x_1}{2}, 0, 0)x$ is of the form ${}^t(0 \ 1 \ \cdots)$. Multiplying an element of $\mathrm{SL}(3)$, we may assume that $x_3 = x_4 = 0$. Since $\Delta(x) \neq 0$, $x_5 \neq 0$. Multiplying an element of $\mathrm{SL}(3)$ again, we may assume $x_6 = x_7 = 0$.

Suppose $x_2 = \cdots = x_7 = 0$. Then the second entry of $u_1(1, 0, 0)x$ is non-zero and it reduces to the previous case. □

**Theorem (4.4)** *The map $\alpha_V : G_k \setminus V_k^{\mathrm{ss}} \to \mathfrak{Er}_2$ is a bijection. If $x \in V_k^{\mathrm{ss}}$ corresponds to $k(\alpha)$ ($\alpha \notin k, \alpha^2 \in k$), $G_x^0 \cong \mathrm{SU}_\alpha(2,1)$. Moreover, for any $x \in V_k^{\mathrm{ss}}$, the corresponding field is $k(\Delta(x)^{\frac{1}{2}})$.*

*Proof.* By the above lemma, any orbit in $G_k \setminus V_k^{\mathrm{ss}}$ is represented by a point of the form $(0 \ 1 \ 0 \ 0 \ x_5 \ 0 \ 0)$, where $x_5 \in k^\times$. It is easy to see that this element belongs to the orbit $G_k w_\alpha$ for some $\alpha$ such that $\alpha^2 \in k$. This implies that the orbits of $w_\alpha$'s which we constructed in §2 exhausts all the orbits in $G_k \setminus V_k^{ss}$. Therefore, $\alpha_V$ is bijective.

Finally, the character $\chi$ is a square of the character $\chi'$. So if $x = gw_\alpha$ and $g \in G_k$, the field corresponding to $x$ is $k(\alpha)$ and $\Delta(x) = \chi'(g)^2 \Delta(w_\alpha)$. This proves $k(\Delta(x)^{\frac{1}{2}}) = k(\Delta(w_\alpha)^{\frac{1}{2}}) = k(\alpha)$. □

This theorem itself was easy to verify and we didn't have to consider Galois cohomology. However, this has the following important interpretation in terms of Galois cohomology, which we need to deal with more interesting case (2).

**Corollary (4.5)** *The map* $\mathrm{H}^1(k, \mathrm{SU}_\alpha(2,1)) \to \mathrm{H}^1(k, G)$ *is injective.*

Next we consider case (2). So we use the definitions of $G, V, w, d(A), \Lambda, \tau, g_\alpha, w_\alpha$ in §3 in (4.6)–(4.16). Note that $g_\alpha \tau g_\alpha^{-1} \in G_{w_\alpha k}$. Therefore, we have a split exact sequence
$$1 \to G^0_{w_\alpha} \to G_{w_\alpha} \to \mathbb{Z}/2\mathbb{Z} \to 1.$$
Since $\mathrm{H}^1(k, \widetilde{T})$ is trivial, this implies
$$\alpha_V^{-1}(k(\alpha)) \cong (\mathbb{Z}/2\mathbb{Z}) \setminus \mathrm{Ker}(\mathrm{H}^1(k, \mathrm{U}_\alpha(1,1)) \to \mathrm{H}^1(k, G)).$$

Let $Y_k$ be the set of non-degenerate $2 \times 2$ Hermitian matrices. As in Proposition (3.3), $\mathrm{H}^1(k, \mathrm{U}_\alpha(1,1))$ is represented by elements of $\mathrm{H}^1(k(\alpha)/k, \mathrm{U}_\alpha(1,1))$. We express elements of $\mathrm{H}^1(k(\alpha)/k, \mathrm{U}_\alpha(1,1))$ as $\{h\}$ where $h \in \mathrm{U}_\alpha(1,1)_{k(\alpha)} \cong \mathrm{GL}(2)_{k(\alpha)}$. Then $\{h\} \to h^{-1}\Lambda$ induces a bijection $\mathrm{H}^1(k, \mathrm{U}_\alpha(1,1)) \cong \mathrm{GL}(2)_{k(\alpha)} \setminus Y_k$. Moreover, the map $\mathrm{H}^1(k, \mathrm{SU}_\alpha(2,1)) \to \mathrm{H}^1(k, \mathrm{U}_\alpha(1,1))$ is identified with
$$Y_k \ni H \to \begin{pmatrix} -(\det H)^{-1} & \\ & H \end{pmatrix} \in X_k.$$

**Lemma (4.6)** *Suppose $H$ is a $2 \times 2$ Hermitian matrix and $\begin{pmatrix} -(\det H)^{-1} & \\ & H \end{pmatrix}$ is $\mathrm{SL}(3)_{k(\alpha)}$-equivalent to $\begin{pmatrix} 1 & \\ & \Lambda \end{pmatrix}$. Then there exists $s \in k^\times$ such that $H$ is $\mathrm{GL}(2)_{k(\alpha)}$-equivalent to $\begin{pmatrix} -s & \\ & 1 \end{pmatrix}$.*

*Proof.* Note by an analogue of the Gram–Schmidt process, we can diagonalize $H$. So we may assume $H = \begin{pmatrix} a & \\ & b \end{pmatrix}$ where $a, b \in k^\times$. Let $\mathrm{N}_{k(\alpha)/k}$ be the norm of $k(\alpha)/k$. By assumption, there exists a non-zero vector $x = {}^t(\begin{array}{ccc} x_1 & x_2 & x_3 \end{array})$ such that
$$-\frac{1}{ab}\mathrm{N}_{k(\alpha)/k}(x_1) + a\mathrm{N}_{k(\alpha)/k}(x_2) + b\mathrm{N}_{k(\alpha)/k}(x_3) = 0.$$

Suppose $x_1 = 0$. Let $z = \frac{x_3}{x_2}$. Then $-\frac{b}{a} = \mathrm{N}_{k(\alpha)/k}(z)$. So
$$\begin{pmatrix} z & \\ & 1 \end{pmatrix} H^* \begin{pmatrix} z & \\ & 1 \end{pmatrix} = \begin{pmatrix} -b & \\ & b \end{pmatrix}.$$

Let
$$A = \begin{pmatrix} 1 & 1 \\ -1 & 1 \end{pmatrix}^{-1} \begin{pmatrix} b^{-1} & \\ & 1 \end{pmatrix} \begin{pmatrix} 1 & 1 \\ -1 & 1 \end{pmatrix}.$$

Then
$$A \begin{pmatrix} -b & \\ & b \end{pmatrix} {}^t A = \begin{pmatrix} -1 & \\ & 1 \end{pmatrix}.$$

Suppose $x_1 \neq 0$. Let $z_1 = \frac{bx_3}{x_1}$, $z_2 = \frac{ax_2}{x_1}$, and $A = \begin{pmatrix} bz_2^\sigma & -az_1^\sigma \\ z_1 & z_2 \end{pmatrix}$. Since $aN_{k(\alpha)/k}(z_1) + bN_{k(\alpha)/k}(z_2) = 1$, $\det A = 1$. It is easy to see $AH^*A = \begin{pmatrix} ab & \\ & 1 \end{pmatrix}$. $\square$

Suppose $\{h\} \in H^1(k, U_\alpha(1,1))$ maps to the trivial element of $H^1(k, G)$. Let $H$ be the corresponding $2 \times 2$ Hermitian matrix. By Corollary (4.5), the Hermitian matrix $\begin{pmatrix} -(\det H)^{-1} & \\ & H \end{pmatrix}$ is $SL(3)_{k(\alpha)}$-equivalent to $\begin{pmatrix} 1 & \\ & \Lambda \end{pmatrix}$. Therefore, by Lemma (4.6), $H$ is equivalent to a matrix of the form $\begin{pmatrix} -s & \\ & 1 \end{pmatrix}$.

Conversely, we show that elements of the above form correspond to orbits in $V_k^{ss}$. For $s \in k^\times$, let $H(s) = \begin{pmatrix} -s & \\ & 1 \end{pmatrix}$ and $\overline{H}(s) = \begin{pmatrix} s^{-1} & & \\ & -s & \\ & & 1 \end{pmatrix}$. Let

(4.7) $$h(s) = \begin{pmatrix} -s^{-1} & \\ & -1 \end{pmatrix}, \quad \overline{h}(s) = \begin{pmatrix} s & & \\ & -s^{-1} & \\ & & -1 \end{pmatrix}.$$

Then $\{h(s)\} \in H^1(k, U_\alpha(2,1))$ and $\{\overline{h}(s)\} \in SU_\alpha(2,1))$ correspond to $H(s)$ and $\overline{H}(s)$ respectively. We define

(4.8) $$w(s) = v_1 \begin{pmatrix} 0 \\ 0 \\ \frac{1+s}{2} \\ \frac{1-s}{2} \\ 0 \\ 0 \\ 0 \end{pmatrix} + v_2 \begin{pmatrix} 0 \\ 0 \\ 0 \\ 0 \\ 0 \\ \frac{1+s}{2} \\ -\frac{1-s}{2} \end{pmatrix}.$$

**Definition (4.9)** $w_\alpha(s) = g_\alpha w(s)$.

**Lemma (4.10)** *The class $\{h(s)\} \in H^1(k, U_\alpha(1,1))$ maps to the trivial element in $H^1(k, G)$ and corresponds to the orbit of $w_\alpha(s)$.*

*Proof.* Let

$$A(s) = \begin{pmatrix} & & 1 \\ & 1 & \\ 1 & & \end{pmatrix} \begin{pmatrix} \frac{1}{2} & -\frac{1}{2} & \\ \frac{1}{2} & \frac{1}{2} & \\ & & 1 \end{pmatrix} \begin{pmatrix} -s^{-1} & & \\ & 1 & \\ & & 1 \end{pmatrix}$$

$$\times \begin{pmatrix} 1 & 1 & \\ -1 & 1 & \\ & & 1 \end{pmatrix} \begin{pmatrix} -s & & \\ & 1 & \\ & & 1 \end{pmatrix}$$

$$= \begin{pmatrix} & & 1 \\ \frac{1+s}{2} & \frac{1-s^{-1}}{2} & \\ \frac{1-s}{2} & \frac{-1-s^{-1}}{2} & \end{pmatrix}.$$

Then $A(s)\overline{H}(s)^*A(s) = \begin{pmatrix} 1 & \\ & \Lambda \end{pmatrix}$. Note that $^*A(s) = {}^tA(s)$. Also

$$^*A(s)^{-1} = \begin{pmatrix} & & 1 \\ \frac{1+s^{-1}}{2} & \frac{1-s}{2} & \\ \frac{1-s^{-1}}{2} & \frac{-1-s}{2} & \end{pmatrix}.$$

Let $B_1(s) = \begin{pmatrix} 1 & \\ & \Lambda \end{pmatrix} {}^*A(s)^{-1} \begin{pmatrix} 1 & \\ & \Lambda \end{pmatrix}$, $B_2(s) = \begin{pmatrix} s & \\ & 1 \end{pmatrix}$, and

(4.11) $$B(s) = \left( \begin{pmatrix} 1 & & \\ & B_1(s) & \\ & & {}^*B_1(s)^{-1} \end{pmatrix}, B_2(s) \right).$$

Then

$$B_1(s)^{-1} \begin{pmatrix} 1 & \\ & \Lambda \end{pmatrix} {}^*B_1(s)^{-1} \begin{pmatrix} 1 & \\ & \Lambda \end{pmatrix} = \begin{pmatrix} 1 & \\ & \Lambda \end{pmatrix} {}^*A(s) \begin{pmatrix} 1 & \\ & \Lambda \end{pmatrix} A(s)$$

$$= \begin{pmatrix} 1 & \\ & \Lambda \end{pmatrix} \overline{H}(s)^{-1} = \overline{h}(s),$$

$$B_2(s)^{-1} \tau_2 B_2(s) \tau_2 = \begin{pmatrix} s^{-1} & \\ & s \end{pmatrix}.$$

So,

(4.12) $$B(s)^{-1} \tau B(s)^\sigma \tau = B(s)^{-1} \tau B(s) \tau = d(h(s)).$$

Therefore, $\{h(s)\}$ maps to the trivial element and it corresponds to the orbit of

$$g_\alpha B(s) g_\alpha^{-1} w_\alpha = g_\alpha B(s) w.$$

Since

$$B_1(s) = \begin{pmatrix} \frac{1+s^{-1}}{2} & \frac{1-s}{2} \\ -\frac{1-s^{-1}}{2} & \frac{1+s}{2} \end{pmatrix}^{-1}, \quad {}^*B_1(s)^{-1} = \begin{pmatrix} \frac{1+s}{2} & \frac{1-s^{-1}}{2} \\ -\frac{1-s}{2} & \frac{1+s^{-1}}{2} \end{pmatrix}^{-1},$$

it is easy to see that $B(s)w = w(s)$. This proves the proposition. □

Now we are ready to determine the fiber structure of $\alpha_V$ for case (2).

**Theorem (4.13)** (1) $\alpha_V^{-1}(k(\alpha)) \cong (\mathbb{Z}/2\mathbb{Z}) \setminus k^\times / \mathrm{N}_{k(\alpha)/k}(k(\alpha)^\times)$.
(2) If $s \in k^\times$, the corresponding orbit is $G_k w_\alpha(s)$, and the non-trivial element of $\mathbb{Z}/2\mathbb{Z}$ acts by $w_\alpha(s) \to w_\alpha(s^{-1})$.
(3) If $x \in \alpha_V^{-1}(k(\alpha))$, $k(\alpha)$ is generated by roots of $F_x(v)$.

*Proof.* We consider $\mathrm{GL}(1)_{k(\alpha)}$ as a group over $k$. We can regard the norm

$$\mathrm{N}_{k(\alpha)/k} : \mathrm{GL}(1)_{k(\alpha)} \to \mathrm{GL}(1)_k$$

as a $k$–homomorphism between $k$–groups. Let $T_\alpha = \mathrm{Ker}(\mathrm{N}_{k(\alpha)/k})$. We can regard $T_\alpha$ as the unitary group of dimension one. So as in the case of $\mathrm{SU}_\alpha(2,1), \mathrm{U}_\alpha(1,1)$, $\mathrm{H}^1(k, T_\alpha)$ is represented by elements of $\mathrm{H}^1(k(\alpha)/k, T_\alpha)$. We can identify $T_{\alpha k(\alpha)} \cong k(\alpha)^\times$, and so we express cohomology classes as $\{s\}$ where $s \in k(\alpha)^\times$. Since the action of $\sigma$ induced from the $k$–structure of $T_\alpha$ is $x \to (x^\sigma)^{-1}$, the cocycle condition is equivalent to $s \in k^\times$ and $s_1, s_2 \in k^\times$ are equivalent if and only if $s_1 s_2^{-1} \in \mathrm{N}_{k(\alpha)/k}(k(\alpha)^\times)$. Therefore, $\mathrm{H}^1(k, T_\alpha) \cong k^\times / \mathrm{N}_{k(\alpha)/k}(k(\alpha)^\times)$. Moreover, the map $\mathrm{H}^1(k, U_\alpha(1,1)) \to \mathrm{H}^1(k, T_\alpha)$ induced by the determinant can be identified with $\{h\} \to \{\det h\}$ for $h \in \mathrm{GL}(2)_{k(\alpha)} \cong \mathrm{U}_\alpha(1,1)_{k(\alpha)}$.

It is easy to see that if $s_1, s_2 \in k^\times$ and $s_1 s_2^{-1} \in \mathrm{N}_{k(\alpha)/k}(k(\alpha)^\times)$, $H(s_1), H(s_2)$ are $\mathrm{GL}(2)_{k(\alpha)}$–equivalent. This implies that $\{h(s_1)\} = \{h(s_2)\}$ in $\mathrm{H}^1(k, \mathrm{U}_\alpha(1,1))$. So the map

$$f : k^\times / \mathrm{N}_{k(\alpha)/k}(k(\alpha)^\times) \ni s \to \{h(s)\} \in \mathrm{Ker}(\mathrm{H}^1(k, \mathrm{U}_\alpha(1,1)) \to \mathrm{H}^1(k, G))$$

is well defined. By Lemma (4.6), $f$ is surjective. Since the composition of $f$ and $\mathrm{H}^1(k, \mathrm{U}_\alpha(1,1)) \to \mathrm{H}^1(k, T_\alpha) \cong k^\times / \mathrm{N}_{k(\alpha)/k}(k(\alpha)^\times)$ is $s \to s^{-1}$, $f$ is bijective.

Since $\tau d(h(s)) \tau = d(h(s^{-1}))$, $\mathbb{Z}/2\mathbb{Z}$ acts on $\mathrm{Ker}(\mathrm{H}^1(k, \mathrm{U}_\alpha(1,1)) \to \mathrm{H}^1(k, G)) \cong k^\times / \mathrm{N}_{k(\alpha)/k}(k(\alpha)^\times)$ by $s \to s^{-1}$. This proves (1), (2).

By straightforward computations,

$$F_{w(s)}(v) = -s v_1 v_2.$$

So the roots of $F_{w_\alpha(s)}(v)$ are $\pm \alpha$. It is easy to see that if $g \in G_k$, the roots of $F_{g w_\alpha(s)}(v)$ and $F_{w_\alpha(s)}(v)$ generate the same field. This proves (3). □

Finally we determine $G^0_{w_\alpha(s)}$.

**Definition (4.14)** $\mathrm{U}_\alpha(H(s)) = \{A \in \mathrm{GL}(2)_{k(\alpha)} \mid A H(s)^* A = \Lambda\}$.

We consider $\mathrm{U}_\alpha(H(s))$ as a group over $k$.

**Proposition (4.15)** $G^0_{w_\alpha(s)} \cong \mathrm{U}_\alpha(H(s)) \times \widetilde{T}$.

*Proof.* Let $R$ be any $k$–algebra. We define $R(\alpha)$ and the action of $\sigma$ on $R(\alpha)$ as in Proposition (2.14).

Any element in $G^0_{w_\alpha(s) R(\alpha)}$ can be expressed as the product of an element of $\widetilde{T}_R$ and

$$g_\alpha B(s) d(g) (g_\alpha B(s))^{-1},$$

where $g \in \mathrm{GL}(2)_{R(\alpha)}$.

Since $g_\alpha^\sigma = g_\alpha \tau$,

$$\begin{aligned}
(g_\alpha B(s)d(g)(g_\alpha B(s))^{-1})^\sigma &= g_\alpha^\sigma B(s)d(g^\sigma)(g_\alpha^\sigma B(s))^{-1} \\
&= g_\alpha \tau B(s)d(g)B(s)^{-1}\tau g_\alpha^{-1} \\
&= g_\alpha B(s)B(s)^{-1}\tau B(s)d(g)B(s)^{-1}\tau B(s)(g_\alpha B(s))^{-1}.
\end{aligned}$$

By straightforward computations using (4.12),

$$(4.16) \qquad B(s)^{-1}\tau B(s) = \left( \begin{pmatrix} 1 & & \overline{H}(s)^{-1} \\ & \overline{H}(s) & \end{pmatrix}, \begin{pmatrix} & -s^{-1} \\ -s & \end{pmatrix} \right).$$

So,
$$B(s)^{-1}\tau B(s)d(g^\sigma)B(s)^{-1}\tau B(s) = d(H(s)^{-1*}g^{-1}H(s)).$$

This implies that

$$(g_\alpha B(s)d(g)(g_\alpha B(s))^{-1})^\sigma = g_\alpha B(s)d(H(s)^{-1*}g^{-1}H(s)))(g_\alpha B(s))^{-1}.$$

Therefore, this element belongs to $G^0_{w_\alpha(s)R}$ if and only if $g = H(s)^{-1*}g^{-1}H(s)$. This condition is equivalent to $g \in \mathrm{U}_\alpha(H(s))_R$. This proves the proposition again by Theorem [7, p. 17]. $\square$

Akihiko Yukie
Oklahoma State University
Mathematics Department
401 Math Science
Stillwater OK 74078-1058 USA
yukie@math.okstate.edu
http://www.math.okstate.edu/~yukie